\documentclass[10pt]{article}

\usepackage{amssymb}
\usepackage{amsmath}
\usepackage{graphicx}
\usepackage{theorem}

\newtheorem{thm}{Theorem}[section]
\newtheorem{lem}{Lemma}[section]
\newtheorem{cor}{Corollary}[section]
\newtheorem{prp}{Proposition}[section]
\theorembodyfont{\rmfamily}
\newtheorem{exm}{Example}[section]
\newtheorem{remark}{Remark}[section]

\makeatletter

\@addtoreset{equation}{section}
\makeatother

\def\ph{{\widehat p}}
\def\Re{{\mathbb{R}}}

\def\Bc{{\cal B}}
\def\Dc{{\cal D}}
\def\Nc{{\cal N}}
\def\Rc{{\cal R}}
\def\Sc{{\cal S}}

\def\al{{\alpha}}
\def\be{{\beta}}
\def\ga{{\gamma}}
\def\de{{\delta}}
\def\la{{\lambda}}
\def\th{{\theta}}
\def\si{{\sigma}}

\def\Ga{{\Gamma}}
\def\De{{\Delta}}
\def\La{{\Lambda}}
\def\Th{{\Theta}}
\def\Si{{\Sigma}}
\def\Om{{\Omega}}
\def\Thh{{\widehat \Th}}

\def\tr{{\rm tr\,}}
\def\diag{{\rm diag\,}}
\def\vec{{\rm vec}}
\def\Er{{\rm E}}
\def\dd{{\rm d}}
\def\infi{{\infty}}

\def\non{{\nonumber}}
\def\qed{\hfill$\Box$}

\begin{document}

\title{Proper Bayes and Minimax Predictive Densities for a Matrix-variate Normal Distribution}
\author{
Hisayuki Tsukuma\footnote{Faculty of Medicine, Toho University, 
5-21-16 Omori-nishi, Ota-ku, Tokyo 143-8540, Japan, E-Mail: tsukuma@med.toho-u.ac.jp}
\ and 
Tatsuya Kubokawa\footnote{Faculty of Economics, University of Tokyo, 
7-3-1 Hongo, Bunkyo-ku, Tokyo 113-0033, Japan, E-Mail: tatsuya@e.u-tokyo.ac.jp}
}
\maketitle 

\begin{abstract}
This paper deals with the problem of estimating predictive densities of a matrix-variate normal distribution with known covariance matrix.
Our main aim is to establish some Bayesian predictive densities related to matricial shrinkage estimators of the normal mean matrix.
The Kullback-Leibler loss is used for evaluating decision-theoretical optimality of predictive densities.
It is shown that a proper hierarchical prior yields an admissible and minimax predictive density.
Also, superharmonicity of prior densities is paid attention to for finding out a minimax predictive density with good numerical performance.

\par\vspace{4mm}\noindent
{\it AMS 2010 subject classifications:} Primary 62C15, 62C20; secondary 62C10.

\par\vspace{2mm}\noindent
{\it Key words and phrases:} Admissibility, Gauss' divergence theorem, generalized Bayes estimator, inadmissibility, Kullback-Leibler loss, minimaxity, shrinkage estimator, statistical decision theory.
\end{abstract}

\section{Introduction}
\label{sec:intro}

The problem of predicting a density function for future observation is an important field in practical applications of statistical methodology.
Since predictive density estimation has been revealed to be parallel to shrinkage estimation for location parameter, it has extensively been studied in the literature.
Particularly, the Bayesian prediction for a multivariate (vector-valued) normal distribution has been developed by Komaki (2001), George et al. (2006) and Brown et al. (2008).
See George et al. (2012) for a broad survey including a clear explanation of parallelism between density prediction and shrinkage estimation.

\medskip
This paper addresses Bayesian predictive density estimation for a matrix-variate normal distribution.
Denote by $\Nc_{a\times b}(M,\Psi\otimes\Si)$ the $a\times b$ matrix-variate normal distribution with mean matrix $M$ and positive definite covariance matrix $\Psi\otimes\Si$, where $M$, $\Psi$ and $\Si$ are, respectively, $a\times b$, $a\times a$ and $b\times b$ matrices of parameters and $\Psi\otimes\Si$ represents the Kronecker product of the positive definite matrices $\Psi$ and $\Si$.
Let $A^\top$ be the transpose of a matrix $A$ and let $\tr A$ and $|A|$ be, respectively, the trace and the determinant a square matrix $A$.
Also, let $A^{-1}$ be the inverse of a nonsingular matrix $A$.
If an $a\times b$ random matrix $Z$ is distributed as $\Nc_{a\times b}(M,\Psi\otimes\Si)$, then $Z$ has density of the form
$$
(2\pi)^{-ab/2}|\Psi|^{-b/2}|\Si|^{-a/2}\exp[-2^{-1}\tr\{\Psi^{-1}(Z-\Th)\Si^{-1}(Z-\Th)^\top\}].
$$
For more details of matrix-variate normal distribution, see Muirhead (1982) and Gupta and Nagar (1999).

\medskip
It is assumed in this paper that the covariance matrix of a matrix-variate normal distribution is known.
Then the prediction problem is more precisely formulated as follows:
Let $X|\Th\sim\Nc_{r\times q}(\Th, v_xI_r\otimes I_q)$ and $Y|\Th\sim\Nc_{r\times q}(\Th, v_yI_r\otimes I_q)$, where $\Th$ is a common $r\times q$ matrix of unknown parameters, $v_x$ and $v_y$ are known positive values and $I_r$ stands for the identity matrix of order $r$.
Assume that $q\geq r$ and $X$ and $Y$ are independent.
Let $p(X\mid \Th)$ and $p(Y\mid \Th)$ be the densities of $X$ and $Y$, respectively.
Consider here the problem of estimating $p(Y\mid \Th)$ based only on the observed $X$.
Denote by $\ph=\ph(Y\mid X)$ an estimated density for $p(Y\mid \Th)$ and hereinafter $\ph$ is referred to as a predictive density of $Y$.
Define the Kullback-Leibler (KL) loss as
\begin{align}\label{eqn:loss}
L_{KL}(\Th,\ph)
&=\Er^{Y|\Th}\bigg[\log {p(Y\mid\Th)\over\ph(Y\mid X)}\bigg] \non\\
&=\int_{\Re^{r\times q}} p(Y\mid \Th)\log {p(Y\mid\Th)\over\ph(Y\mid X)}\dd Y.
\end{align}
The performance of a predictive density $\ph$ is evaluated by the risk function with respect to the KL loss (\ref{eqn:loss}),
\begin{align*}
R_{KL}(\Th,\ph)&=\Er^{X|\Th}[L_{KL}(\Th,\ph)]\\
&=\int_{\Re^{r\times q}}\int_{\Re^{r\times q}}p(X\mid\Th)p(Y\mid\Th)\log {p(Y\mid\Th)\over\ph(Y\mid X)}\dd Y\dd X.
\end{align*}

\medskip
Let $\pi(\Th)$ be a proper/improper density of prior distribution for $\Th$, where we assume that the marginal density of $X$,
$$
m_\pi(X;v_x)=\int_{\Re^{r\times q}} p(X\mid \Th)\pi(\Th)\dd\Th,
$$
is finite for all $X\in\Re^{r\times q}$.
Denote the Frobenius norm of a matrix $A$ by $\Vert A\Vert=\sqrt{\tr AA^\top}$.
Let
$$
p_\pi(X,Y)=\int_{\Re^{r\times q}} p(X\mid \Th) p(Y\mid \Th)\pi(\Th)\dd\Th.
$$
Note that $p_\pi(X,Y)$ is finite if $m_\pi(X;v_x)$ is finite.
Here $p_\pi(X,Y)$ can be rewritten as
\begin{align*}
p_\pi(X,Y)&={1\over (2\pi v_s)^{qr/2}}e^{-\Vert Y-X\Vert^2/2v_s} \times
\int_{\Re^{r\times q}} {1\over (2\pi v_w)^{qr/2}}e^{-\Vert W-\Th\Vert^2/2v_w}\pi(\Th)\dd\Th \\
&\equiv \ph_U(Y\mid X) \times  m_\pi(W;v_w),
\end{align*}
where $v_s=v_x+v_y$ and
$$
W=v_w(X/v_x+Y/v_y)\mid \Th \sim\Nc_{r\times q}(\Th, v_wI_r\otimes I_q)
$$
with $v_w=(1/v_x+1/v_y)^{-1}$.
From Aitchison (1975), a Bayesian predictive density relative to the KL loss (\ref{eqn:loss}) is given by
\begin{equation}\label{eqn:BPD}
\ph_\pi(Y\mid X)={p_\pi(X,Y) \over m_\pi(X; v_x)}
={m_\pi(W;v_w)\over m_\pi(X;v_x)}\,\ph_U(Y\mid X).
\end{equation}
See George et al. (2006, Lemma 2) for the multivariate (vector-valued) normal case.

\medskip
It is noted that $\ph_U(Y\mid X)$ is the Bayesian predictive density with respect to the uniform prior $\pi_U(\Th)=1$.
Under the predictive density estimation problem relative to the KL loss (\ref{eqn:loss}), $\ph_U(Y\mid X)$ is the best invariant predictive density with respect to a location group.
Using the same arguments as in George et al. (2006, Corollary 1) gives that, for any $r$ and $q$, $\ph_U(Y\mid X)$ is minimax relative to the KL loss (\ref{eqn:loss}) and has a constant risk.

\medskip
Recently, Matsuda and Komaki (2015) constructed an improved Bayesian predictive density on $\ph_U(Y\mid X)$ by using a prior density of the form
\begin{equation}\label{eqn:pr_em}
\pi_{EM}(\Th)=|\Th\Th^\top|^{-\al^{EM}/2},\quad \al^{EM}=q-r-1.
\end{equation}
The prior (\ref{eqn:pr_em}) is interpreted as an extension of Stein's (1973, 1981) harmonic prior
\begin{equation}\label{eqn:pr_js}
\pi_{JS}(\Th)=\Vert\Th\Vert^{-\be^{JS}}=\{\tr(\Th\Th^\top)\}^{-\be^{JS}/2},\quad \be^{JS}=qr-2.
\end{equation}
In the context of Bayesian estimation for mean matrix, (\ref{eqn:pr_em}) yields a matricial shrinkage estimator, while (\ref{eqn:pr_js}) does a scalar shrinkage one.
Note that, when $X\sim\Nc_{r\times q}(\Th, v_x I_r\otimes I_q)$, typical examples of the matricial and the scalar shrinkage estimators for $\Th$ are, respectively, the Efron-Morris (1972) estimator 
\begin{equation}\label{eqn:EM}
\Thh_{EM}=\{I_r-\al^{EM}v_x(XX^\top)^{-1}\}X
\quad \textup{for $\al^{EM}\geq 1$ (i.e., $q\geq r+2$)}
\end{equation}
and the James-Stein (1961) like estimator
\begin{equation}\label{eqn:JS}
\Thh_{JS}=\Big\{1-\frac{\be^{JS}v_x}{\tr(XX^\top)}\Big\}X
\quad \textup{for $\be^{JS}\geq 1$ (i.e., $qr\geq 3$)}.
\end{equation}
The two estimators $\Thh_{EM}$ and $\Thh_{JS}$ are minimax relative to a quadratic loss.
Also, $\Thh_{EM}$ and $\Thh_{JS}$ are characterized as empirical Bayes estimators, but they are not generalized Bayes estimators which minimize the posterior expected quadratic loss.

\medskip
The purposes of this paper are to construct some Bayesian predictive densities with different priors from (\ref{eqn:pr_em}) and (\ref{eqn:pr_js}) and to discuss their decision-theoretic properties such as admissibility and minimaxity.
Section \ref{sec:preliminaries} first lists some results on the Kullback-Leibler risk and the differentiation operators.
Section \ref{sec:properminimax} applies an extended Faith's (1978) prior to our predictive density estimation problem and provides sufficient conditions for minimaxity of the resulting Bayesian predictive densities.
Also, an admissible and minimax predictive density is obtained by considering a proper hierarchical prior. 
In Section \ref{sec:superharmonic}, we utilize Stein's (1973, 1981) ideas for deriving some minimax predictive densities with superharmonic priors.
Section \ref{sec:MCstudies} investigates numerical performance in risk of some Bayesian minimax predictive densities.

\section{Preliminaries}\label{sec:preliminaries}
\subsection{The Kullback-Leibler risk}

First, we state some useful lemmas in terms of the Kullback-Leibler (KL) risk.
The lemmas are based on Stein (1973, 1981), George et al. (2006) and Brown et al. (2008) and play important roles in studying decision-theoretic properties of a Bayesian predictive density.

\medskip
From George et al. (2006, Lemma 3), we observe that $m_\pi(W;v_w)<\infi$ for all $W\in\Re^{r\times q}$ if $m_\pi(X;v_x)<\infi$ for all $X\in\Re^{r\times q}$.
Note also that $\int_{\Re^{r\times q}}\ph_\pi(Y\mid X)\dd Y=1$ and
$$
\int_{\Re^{r\times q}}Y\ph_\pi(Y\mid X)\dd Y
=\frac{\int_{\Re^{r\times q}} \Th p(X\mid \Th)\pi(\Th)\dd\Th}{\int_{\Re^{r\times q}} p(X\mid \Th)\pi(\Th)\dd\Th},
$$
namely, the mean of a predictive distribution for $Y$ is the same as the posterior mean of $\Th$ given $X$ or, equivalently, the generalized Bayes estimator relative to a quadratic loss for a mean of $X$.

\medskip
Hereafter denote by $p(W|\Th)$ a density of $W|\Th\sim\Nc_{r\times q}(\Th, vI_r\otimes I_q)$ with a positive value $v$.
In order to prove minimaxity of a Bayesian predictive density, we require the following lemma, which implies that our Bayesian prediction problem can be reduced to the Bayesian estimation problem for the normal mean matrix relative to a quadratic loss.
\begin{lem}\label{lem:identity}
The KL risk difference between $\ph_U(Y\mid X)$ and $\ph_\pi(Y\mid X)$ can be written as
$$
R_{KL}(\Th,\ph_U) - R_{KL}(\Th,\ph_\pi) 
=\frac{1}{2}\int_{v_w}^{v_x}\frac{1}{v^2}\{\Er^{W|\Th}[\Vert W-\Th\Vert^2]-\Er^{W|\Th}[\Vert\Thh_\pi-\Th\Vert^2]\}\dd v,
$$
where $\Er^{W|\Th}$ stands for expectation with respect to $W$ and
$$
\Thh_\pi=\Thh_\pi(W)=\frac{\int_{\Re^{r\times q}} \Th p(W\mid \Th)\pi(\Th)\dd\Th}{\int_{\Re^{r\times q}} p(W\mid \Th)\pi(\Th)\dd\Th}.
$$
\end{lem}

{\bf Proof.}\ \ 
This is verified by the same arguments as in Brown et al. (2008, Theorem 1 and its proof).
\hfill$\Box$

\medskip
Let $\nabla_W=(\partial/\partial w_{ij})$ be an $r\times q$ matrix of differentiation operators with respect to an $r\times q$ matrix $W=(w_{ij})$ of full row rank.
For a scalar function $g(W)$ of $W$, the operation $\nabla_W g(W)$ is defined as an $r\times q$ matrix whose $(i,j)$-th element is $\partial g(W)/\partial w_{ij}$.
Also, for a $q\times a$ matrix-valued function $G(W)=(g_{ij})$ of $W$, the operation $\nabla_WG(W)$ are defined as an $r\times a$ matrix whose $(i,j)$-th element of $\nabla_WG(W)$ is $\sum_{k=1}^q \partial g_{kj}/\partial w_{ik}$.

\medskip
Stein (1973) showed that for a $q\times r$ matrix $G(W)$
$$
\int_{\Re^{r\times q}} \tr\{(W-\Th)G(W)\}p(W|\Th)\dd W = v \int_{\Re^{r\times q}}\tr\{\nabla_WG(W)\}p(W|\Th)\dd W,
$$
namely,
$$
\Er^{W|\Th}[\tr\{(W-\Th)G(W)\}] = v \Er^{W|\Th}[\tr\{\nabla_WG(W)\}].
$$
This identity is referred to as the Stein identity in the literature.
Using the Stein identity, we can easily obtain the following lemma.
\begin{lem}\label{lem:identity2}
Use the same notation as in Lemma \ref{lem:identity}.
Then we obtain
\begin{align*}
&R_{KL}(\Th,\ph_U) - R_{KL}(\Th,\ph_\pi)\\
&=-\int_{v_w}^{v_x}\Er^{W|\Th}\bigg[2\frac{\tr[\nabla_W\nabla_W^\top m_\pi(W;v)]}{m_\pi(W;v)}-\frac{\Vert\nabla_W m_\pi(W;v)\Vert^2}{\{m_\pi(W;v)\}^2}\bigg]\dd v.
\end{align*}
\end{lem}

{\bf Proof.}\ \ 
This lemma can be shown by the same arguments as in Stein (1973, 1981).
We provide only an outline of proof.

\medskip
Note from Brown (1971) that $\Thh_\pi$, given in Lemma \ref{lem:identity}, can be represented as
$$
\Thh_\pi=W+v\frac{\nabla_W m_\pi(W;v)}{m_\pi(W;v)}=W+v\nabla_W\log m_\pi(W;v).
$$
By some manipulation after using the Stein identity, we have
\begin{align*}
&\frac{1}{v}\{\Er^{W|\Th}[\Vert W-\Th\Vert^2]-\Er^{W|\Th}[\Vert\Thh_\pi-\Th\Vert^2]\} \\
&=-v\Er^{W|\Th}\bigg[2\frac{\tr[\nabla_W\nabla_W^\top m_\pi(W;v)]}{m_\pi(W;v)}-\frac{\Vert\nabla_W m_\pi(W;v)\Vert^2]}{\{m_\pi(W;v)\}^2}\bigg].
\end{align*}
Combining this identity and Lemma \ref{lem:identity} completes the proof.
\hfill$\Box$

\medskip
Using Lemma \ref{lem:identity2} immediately establishes the following proposition.
\begin{prp}\label{prp:cond_mini}
$\ph_\pi(Y|X)$ is minimax relative to the KL loss (\ref{eqn:loss}) if
$$
2\tr[\nabla_W\nabla_W^\top m_\pi(W;v)]-\frac{\Vert\nabla_W m_\pi(W;v)\Vert^2}{m_\pi(W;v)}\leq 0
$$
for $v_w\leq v\leq v_x$.
\end{prp}

\subsection{Differentiation of matrix-valued functions}

Next, some useful formulae are listed for differentiation with respect to a symmetric matrix.
The formulae are applied to evaluation of the Kullback-Leibler risks of our Bayesian predictive densities.

\medskip
Let $S=(s_{ij})$ be an $r\times r$ symmetric matrix of full rank.
Let $\Dc_S$ be an $r\times r$ symmetric matrix of differentiation operators with respect to $S$, where the $(i,j)$-th  element of $\Dc_S$ is
$$
\{\Dc_S\}_{ij}=\frac{1+\de_{ij}}{2}\frac{\partial}{\partial s_{ij}}
$$
with the Kronecker delta $\de_{ij}$.

\medskip
Let $g(S)$ be a scalar-valued and differentiable function of $S=(s_{ij})$.
Also let $G(S)=(g_{ij}(S))$ be an $r\times r$ matrix, where all the elements $g_{ij}(S)$ are differentiable functions of $S$.
The operations $\Dc_S g(S)$ and $\Dc_S G(S)$ are, respectively, $r\times r$ matrices, where the $(i,j)$-th elements of $\Dc_S g(S)$ and $\Dc_S G(S)$ are defined as, respectively,
$$
\{\Dc_S g(S)\}_{ij}=\frac{1+\de_{ij}}{2}\frac{\partial g(S)}{\partial s_{ij}},\quad 
\{\Dc_S G(S)\}_{ij}=\sum_{k=1}^r\frac{1+\de_{ik}}{2}\frac{\partial g_{kj}(S)}{\partial s_{ik}}.
$$

\medskip
First, the product rule in terms of $\Dc_S$ is expressed in the following lemma due to Haff (1982).
\begin{lem}\label{lem:diff1}
Let $G_1$ and $G_2$ be $r\times r$ matrices such that all the elements of $G_1$ and $G_2$ are differentiable functions of $S$.
Then we have
$$
\Dc_S (G_1G_2)=(\Dc_S G_1)G_2+(G_1^\top \Dc_S)^\top G_2.
$$ 
In particular, for differentiable scalar-valued functions $g_1(S)$ and $g_2(S)$,
$$
\Dc_S \{g_1(S)g_2(S)\}=g_2(S)\Dc_S g_1(S)+g_1(S) \Dc_S g_2(S).
$$ 
\end{lem}

\medskip
Denote by $S=HLH^\top$ the eigenvalue decomposition of $S$, where $H=(h_{ij})$ is an orthogonal matrix of order $r$ and $L=\diag(\ell_1,\ldots,\ell_r)$ is a diagonal matrix of order $r$ with $\ell_1\geq \cdots\geq \ell_r$.
The following lemma is provided by Stein (1973).
\begin{lem}\label{lem:diff2}
Define $\Psi(L)=\diag(\psi_1,\ldots,\psi_r)$, whose diagonal elements are differentiable functions of $L$.
Then we obtain 
\begin{enumerate}
\parskip=0pt
\item[{\rm (i)}] $\{\Dc_S\}_{ij} \ell_k=h_{ik}h_{jk}$\ \ $(k=1,\ldots,r)$,
\item[{\rm (ii)}] $\Dc_S H\Psi(L)H^t=H\Psi^*(L)H^t$, where $\Psi^*(L)=\diag(\psi_1^*,\ldots,\psi_r^*)$ with
$$
\psi_i^*=\frac{\partial \psi_i}{\partial\ell_i}+\frac{1}{2}\sum_{j\ne i}^r\frac{\psi_i-\psi_j}{\ell_i-\ell_j}.
$$
\end{enumerate}
\end{lem}

\begin{lem}\label{lem:diff3}
Let $a$ and $b$ be constants and let $C$ be a symmetric constant matrix $C$.
Then it holds that
\begin{enumerate}
\parskip=0pt
\item[{\rm (i)}] $\Dc_S \tr(S C)=C$,
\item[{\rm (ii)}] $\displaystyle \Dc_S S=\frac{r+1}{2}I_r$,
\item[{\rm (iii)}] $\displaystyle \Dc_S S^2=\frac{r+2}{2}S+\frac{1}{2}(\tr S)I_r$.
\item[{\rm (iv)}] $\Dc_S |aI_r+bS|=b|aI_r+bS|(aI_r+bS)^{-1}$ if $aI_r+bS$ is nonsingular.
\end{enumerate}
\end{lem}

{\bf Proof.}\ \ 
For proofs of Parts (i), (ii) and (iii), see Haff (1982) and Magnus and Neudecker (1999).
Using (i) of Lemma \ref{lem:diff2} gives that
\begin{align*}
\{\Dc_S |aI_r+bS|\}_{ij} &= \{\Dc_S\}_{ij} \prod_{k=1}^r(a+b\ell_k) \\
&=b\sum_{c=1}^r h_{ic}h_{jc}\prod_{k\ne c}^r(a+b\ell_k) \\
&=b|aI_r+bS|\sum_{c=1}^r h_{ic}h_{jc}(a+b\ell_c)^{-1} \\
&=b|aI_r+bS|\{(aI_r+bS)^{-1}\}_{ij},
\end{align*}
which implies Part (iv).
\qed

\medskip
Let $\nabla_W$ be the same $r\times q$ differentiation operator matrix as in the preceding subsection.
If $S=WW^\top$, then we have the following lemma, where the proof is referred to in Konno (1992).
\begin{lem}\label{lem:diff4}
Let $G$ be an $r\times r$ symmetric matrix, where all the elements of $G$ are differentiable function of $S=WW^\top$.
Then it holds that
\begin{enumerate}
\parskip=0pt
\item[{\rm (i)}] $\nabla_W^\top G=2W^\top \Dc_S G$,
\item[{\rm (ii)}] $\tr(\nabla_W W^\top G)=(q-r-1)\tr G+2\tr(\Dc_S SG)$.
\end{enumerate}
\end{lem}

\section{Admissible and minimax predictive densities}
\label{sec:properminimax}

In this section, we consider a class of hierarchical priors inspired by Faith (1978) and derive a sufficient condition for minimaxity of the resulting Bayesian predictive density.
Also, a proper Bayes and minimax predictive density is provided.

\subsection{A class of hierarchical prior distributions}

Let $\Sc_r$ be the set of $r\times r$ symmetric matrices.
For $A$ and $B\in\Sc_r$, write $A\prec(\preceq) B$ or $B\succ(\succeq) A$ if $B-A$ is a positive (semi-)definite matrix.
The set $\Rc_r$ is defined as
$$
\Rc_r=\{ \La\in\Sc_r \mid 0_{r\times r}\prec \La \prec I_r\},
$$
where $0_{r\times r}$ is the $r\times r$ zero matrix.
Denote the boundary of $\Rc_r$ by $\partial\Rc_r$.
It is noted that if $\Om\in\partial\Rc_r$ then $0_{r\times r}\preceq \Om \preceq I_r$ and also then $|\Om|=0$ or $|I_r-\Om|=0$.

\medskip
Consider a proper/improper hierarchical prior
$$
\pi_H(\Th)=\int_{\Rc_r}\pi_1(\Th|\Om)\pi_2(\Om)\dd\Om.
$$
The priors $\pi_1(\Th|\Om)$ and $\pi_2(\Om)$ are specified as follows:
Assume that a prior distribution of $\Th$ given $\Om$ is $\Nc_{r\times q}(0_{r\times q},v_0\Om^{-1}(I_r-\Om)\otimes I_q)$, where $v_0$ is a known constant satisfying 
$$
v_0\geq v_x.
$$
Then the first-stage prior density $\pi_1(\Th|\Om)$ can be written as
\begin{equation}\label{eqn:pr_Th}
\pi_1(\Th|\Om)=(2\pi v_0)^{-qr/2}|\Om(I_r-\Om)^{-1}|^{q/2}\exp\Big[-\frac{1}{2v_0}\tr\{\Om(I_r-\Om)^{-1}\Th\Th^\top\}\Big].
\end{equation}
Assume also that $\pi_2(\Om)$, a second-stage prior density for $\Om$, is a differentiable function on $\Rc_r$.

\medskip
Denote by $\ph_H=\ph_H(Y|X)$ the resulting Bayesian predictive density with respect to the hierarchical prior $\pi_H(\Th)$.
Assume that a marginal density of $W$ with respect to $\pi_H(\Th)$ is finite when $v=v_x$.
The marginal density is given by
\begin{align}\label{eqn:m(W)}
m(W)&=\int_{\Re^{r\times q}} p(W|\Th)\pi_H(\Th)\dd\Th \non\\
&=\int_{\Rc_r}\int_{\Re^{r\times q}} \pi(\Th|\Om,W) \dd\Th \pi_2(\Om)\dd\Om,
\end{align}
where $\pi(\Th|\Om,W)=p(W|\Th)\pi_1(\Th|\Om)$ is a posterior density of $\Th$ given $\Om$ and $W$.
To make it easy to derive sufficient conditions that $\ph_H$ is minimax, we show the following lemma.
\begin{lem}\label{lem:alter_m(W)}
The marginal density $m(W)$ can alternatively be represented as
$$
m(W)=\int_{\Rc_r}f_\pi(\La;W)\dd\La,
$$
where
$$
f_\pi(\La;W)=(2\pi v)^{-qr/2}|\La|^{q/2}\pi_2^J(\La)\exp\Big[-\frac{1}{2v}\tr(\La WW^\top)\Big]
$$
with
$$
\pi_2^J(\La)
=v_1^{r(r+1)/2}|v_1I_r+(1-v_1)\La|^{-r-1}\pi_2[\La\{v_1I_r+(1-v_1)\La\}^{-1}].
$$
\end{lem}

{\bf Proof.}\ \ 
Let
$$
\La(I_r-\La)^{-1}=v_1\Om(I_r-\Om)^{-1},\quad v_1=\frac{v}{v_0},
$$
where $0_{r\times r}\prec\La\prec I_r$.
Since $v^{-1}(I_r-\La)^{-1}=v^{-1}I_r+v_0^{-1}\Om(I_r-\Om)^{-1}$, we observe that
\begin{align*}
&\frac{1}{v}\Vert W-\Th\Vert^2+\frac{1}{v_0}\tr\{\Om(I_r-\Om)^{-1}\Th\Th^\top\}\\
&=\frac{1}{v}\tr\Big[(I_r-\La)^{-1}\{\Th-(I_r-\La)W\}\{\Th-(I_r-\La)W\}^\top\Big]
 +\frac{1}{v}\tr(\La WW^\top),
\end{align*}
so $\pi(\Th|\Om,W)$ is proportional to
$$
\pi(\Th|\Om,W)
\propto \exp\Big[-\frac{1}{2v}\tr\Big[(I_r-\La)^{-1}\{\Th-(I_r-\La)W\}\{\Th-(I_r-\La)W\}^\top\Big]\Big],
$$
namely, $\Th|\Om,W\sim\Nc_{r\times q}((I_r-\La)W,v(I_r-\La)\otimes I_q)$.
Integrating out (\ref{eqn:m(W)}) with respect to $\Th$ gives that
\begin{equation}\label{eqn:m(W)-1}
m(W)=(2\pi v)^{-qr/2}\int_{\Rc_r} |\La|^{q/2}\pi_2(\Om)\exp\Big[-\frac{1}{2v}\tr(\La WW^\top)\Big]\dd\Om.
\end{equation}
Note that $\Om=\La\{v_1I_r+(1-v_1)\La\}^{-1}$ and the Jacobian of the transformation from $\Om$ to $\La$ is given by
$$
J[\Om\to \La]=v_1^{r(r+1)/2}|v_1I_r+(1-v_1)\La)|^{-r-1}.
$$
Hence making the transformation from $\Om$ to $\La$ for (\ref{eqn:m(W)-1}) completes the proof.
\qed

\medskip
Let $\Dc_\La$ be an $r\times r$ symmetric matrix of differentiation operators with respect to $\La=(\la_{ij})$, where the $(i,j)$-th  element of $\Dc_\La$ is
$$
\{\Dc_\La\}_{ij}=\frac{1+\de_{ij}}{2}\frac{\partial}{\partial\la_{ij}}.
$$
Proposition \ref{prp:cond_mini} and Lemma \ref{lem:alter_m(W)} are utilized to get sufficient conditions for minimaxity of $\ph_H$.
\begin{thm}\label{thm:faith}
Let $f_\pi(\La;W)$ and $\pi_2^J(\La)$ be defined as in Lemma \ref{lem:alter_m(W)}.
Let
$$
M=M(W)=\int_{\Rc_r}\La f_\pi(\La;W)\dd\La.
$$
Assume that
$$
f_\pi(\La;W)=0 \quad \textup{for all $\La\in\partial\Rc_r$}.
$$
Then $\ph_H$ is minimax relative to the KL loss $(\ref{eqn:loss})$ if $\De(W;\pi_2^J)\leq 0$, where
$$
\De(W;\pi_2^J)=\De_1(W;\pi_2^J)-\De_2(W;\pi_2^J)-(q-3r-3)\tr M
$$
with
\begin{align*}
\De_1(W;\pi_2^J)&=4\int_{\Rc_r}\frac{1}{\pi_2^J(\La)}\tr\{\La^2\Dc_\La \pi_2^J(\La)\}f_\pi(\La;W)\dd\La,\\
\De_2(W;\pi_2^J)&=\frac{2}{m(W)}\int_{\Rc_r}\frac{1}{\pi_2^J(\La)}\tr\{\La M\Dc_\La \pi_2^J(\La)\}f_\pi(\La;W)\dd\La,
\end{align*}
provided all the integrals are finite.
\end{thm}

{\bf Proof.}\ \ 
From Proposition \ref{prp:cond_mini}, $\ph_H$ is minimax when 
$$
\De=2\tr[\nabla_W\nabla_W^\top m(W)]-\frac{\Vert\nabla_W m(W)\Vert^2}{m(W)}\leq 0.
$$

\medskip
It is seen from Lemma \ref{lem:alter_m(W)} that
$$
\nabla_W f_\pi(\La;W)=-\frac{1}{v}\La W f_\pi(\La;W)
$$
and
$$
\nabla_W\nabla_W^\top f_\pi(\La;W)=\Big(\frac{1}{v^2}\La WW^\top\La-\frac{q}{v}\La\Big)f_\pi(\La;W).
$$
Hence we obtain
\begin{equation}\label{eqn:d2_mw1}
\De=\frac{1}{v}[2E_1(W)-\{m(W)\}^{-1}E_2(W)],
\end{equation}
where
\begin{align*}
E_1(W)&=\int_{\Rc_r}\Big[\frac{1}{v}\tr(WW^\top\La^2)-q\tr\La\Big]f_\pi(\La;W)\dd\La,\\
E_2(W)&=\frac{1}{v}\tr\bigg[WW^\top\bigg\{\int_{\Rc_r}\La f_\pi(\La;W)\dd\La\bigg\}^2\bigg] \\
&=\frac{1}{v}\int_{\Rc_r}\tr(MWW^\top\La)f_\pi(\La;W)\dd\La.
\end{align*}

\medskip
Using Lemmas \ref{lem:diff1} and \ref{lem:diff3} yields that
$$
\Dc_\La f_\pi(\La;W)=\frac{1}{2}\Big[q\La^{-1}+\frac{2}{\pi_2^J(\La)}\Dc_\La \pi_2^J(\La)-\frac{1}{v}WW^\top\Big]f_\pi(\La;W),
$$
so that
\begin{align*}
\tr[\Dc_\La\{f_\pi(\La;W)\La^2\}]&=\tr[\La^2\Dc_\La f_\pi(\La;W)]+f_\pi(\La;W)\tr[\Dc_\La\La^2]\\
&=\frac{1}{2}\Big[\frac{2}{\pi_2^J(\La)}\tr[\La^2\Dc_\La \pi_2^J(\La)]-\Big\{\frac{1}{v}\tr(WW^\top\La^2)-q\tr\La\Big\} \\
&\qquad +2(r+1)\tr\La \Big]f_\pi(\La;W).
\end{align*}
Thus $E_1(W)$ can be expressed as
\begin{align}\label{eqn:E1}
E_1(W)
&=2(r+1)\tr M+\frac{1}{2}\De_1(W;\pi_2^J) -2\int_{\Rc_r}\tr[\Dc_\La\{f_\pi(\La;W)\La^2\}]\dd\La.
\end{align}

Similarly, we observe that from Lemmas \ref{lem:diff1} and \ref{lem:diff3}
\begin{align*}
\tr[\Dc_\La\{f_\pi(\La;W)\La\}M]
&=\tr[\La M\Dc_\La f_\pi(\La;W)]+f_\pi(\La;W)\tr[M\Dc_\La\La]\\
&=\frac{1}{2}\Big[(q+r+1)\tr M+\frac{2}{\pi_2^J(\La)}\tr[\La M\Dc_\La \pi_2^J(\La)] \\
&\qquad -\frac{1}{v}\tr(MWW^\top\La)\Big]f_\pi(\La;W),
\end{align*}
which leads to
\begin{align}\label{eqn:E2}
E_2(W)
&=(q+r+1)m(W)\tr M +m(W)\De_2(W;\pi_2^J) \non\\
&\qquad -2\int_{\Rc_r}\tr[\Dc_\La\{f_\pi(\La;W)\La\}M]\dd\La .
\end{align}
Combining (\ref{eqn:d2_mw1}), (\ref{eqn:E1}) and (\ref{eqn:E2}) gives that
\begin{align}\label{eqn:d2_mw2}
\De 
&= \frac{\De(W;\pi_2^J)}{v} -\frac{4}{v}\int_{\Rc_r}\tr[\Dc_\La\{f_\pi(\La;W)\La^2\}]\dd\La \non\\
&\qquad  +\frac{2}{vm(W)}\int_{\Rc_r}\tr[\Dc_\La\{f_\pi(\La;W)\La\}M]\dd\La.
\end{align}
If we can show that two integrals in (\ref{eqn:d2_mw2}) are, respectively, equal to zero, then the proof is complete.

\medskip
Let $G=(g_{ij})$ be an $r\times r$ symmetric matrix such that all the elements of $G$ are differentiable functions of $\La\in\Rc_r$.
Denote
$$
\vec(G)=(g_{11},g_{12},\ldots,g_{1r},g_{22},g_{23},\ldots,g_{r-1,r-1},g_{r-1,r},g_{rr})^\top,
$$
which is a $\{2^{-1}r(r+1)\}$-dimensional column vector.
Denote an outward unit normal vector at a point $\La$ on $\partial\Rc_r$ by
$$
\nu=\nu(\La)=(\nu_{11},\nu_{12},\ldots,\nu_{1r},\nu_{22},\nu_{23},\ldots,\nu_{r-1,r-1},\nu_{r-1,r},\nu_{rr})^\top.
$$
If $\tr(\Dc_\La G)$ is integrable on $\Rc_r$ then it is seen that
$$
\int_{\Rc_r}\tr(\Dc_\La G)\dd\La
=\int_{\Rc_r}\sum_{i=1}^r\sum_{j=1}^r\frac{1+\de_{ij}}{2}\frac{\partial g_{ji}}{\partial\la_{ij}}\dd\La
=\int_{\Rc_r}\sum_{i=1}^r\sum_{j=i}^r\frac{\partial g_{ij}}{\partial\la_{ij}}\dd\La
$$
by symmetry of $\La$ and $G$.
From the Gauss divergence theorem, we obtain
$$
\int_{\Rc_r}\tr(\Dc_\La G)\dd\La=\int_{\partial\Rc_r}\sum_{i=1}^r\sum_{j=i}^r\nu_{ij} g_{ij}\dd\si=\int_{\partial\Rc_r}\nu^\top\vec(G)\dd\si,
$$
where $\si$ stands for Lebesgue measure on $\partial\Rc_r$.

\medskip
Note that
\begin{align*}
\tr[\Dc_\La\{f_\pi(\La;W)\La\}M]=\tr[\Dc_\La\{f_\pi(\La;W)\La M\}]=\tr[\Dc_\La\{f_\pi(\La;W)M\La\}]
\end{align*}
because $M=M(W)$ is symmetric and does not depend on $\La$.
It is observed that $\La^2$ and $\La M+M\La$ are symmetric for $\La\in\Rc_r$, so that
\begin{equation}\label{eqn:gdt1}
\int_{\Rc_r}\tr[\Dc_\La\{f_\pi(\La;W)\La^2\}]\dd\La=\int_{\partial\Rc_r}\nu^\top\vec(\La^2)f_\pi(\La;W)\dd\si,
\end{equation}
and
\begin{align}
&\int_{\Rc_r}\tr[\Dc_\La\{f_\pi(\La;W)\La\}M]\dd\La \non\\
&=\frac{1}{2}\int_{\Rc_r}\tr[\Dc_\La\{f_\pi(\La;W)\La M+f_\pi(\La;W)M\La\}]\dd\La \non\\
&=\frac{1}{2}\int_{\partial\Rc_r}\nu^\top \vec(\La M+M\La)f_\pi(\La;W)\dd\si. 
\label{eqn:gdt2}
\end{align}
Recall that $M$ is finite and $0_{r\times r} \preceq \La \preceq I_r$ for $\La\in\partial\Rc_r$, so that $\nu^\top\vec(\La^2)$ and $\nu^\top \vec(\La M+M\La)$ are bounded.
Since $f_\pi(\La;W)=0$ for any $\La\in\partial\Rc_r$, (\ref{eqn:gdt1}) and (\ref{eqn:gdt2}) are, respectively, equal to zero, which completes the proof. 
\hfill$\Box$

\subsection{Proper Bayes and minimax predictive densities}

Define a second-stage prior density for $\Om$ as
\begin{equation}\label{eqn:pr_GB}
\pi_{GB}(\Om)=K_{a,b}|\Om|^{a/2-1}|I_r-\Om|^{b/2-1},\qquad 0_{r\times r}\prec \Om\prec I_r,
\end{equation}
where $a$ and $b$ are constants and $K_{a,b}$ is a normalizing constant.
The hierarchical prior (\ref{eqn:pr_Th}) with (\ref{eqn:pr_GB}) is a generalization of Faith (1978) in Bayesian minimax estimation of a normal mean vector.
Faith's (1978) prior has also been discussed in detail by Maruyama (1998).
When $a>0$ and $b>0$, $\pi_{GB}(\Om)$ is proper and the distribution of $\Om$ is often called the matrix-variate beta distribution.
Konno (1988) showed that
\begin{align*}
&\int_{\Rc_r} \Om \pi_{GB}(\Om)\dd\Om=\frac{a+r-1}{a+b+2r-2}I_r \qquad \textup{for $a>0$ and $b>0$},\\
&\int_{\Rc_r} \Om(I_r-\Om)^{-1} \pi_{GB}(\Om)\dd\Om=\frac{a+r-1}{b-2}I_r \qquad \textup{for $a>0$ and $b>2$}.
\end{align*}
For other properties of the matrix-variate beta distribution, see Muirhead (1982) and Gupta and Nagar (1999).

\medskip
Let $\ph_{GB}(Y\mid X)$ be the generalized Bayesian predictive density with respect to (\ref{eqn:pr_Th}) and (\ref{eqn:pr_GB}).
A sufficient condition for minimaxity of $\ph_{GB}(Y\mid X)$ is given as follows.
\begin{prp}\label{prp:mini_GB}
Assume that $q-r-1>0$.
Then $\ph_{GB}(Y\mid X)$ is minimax relative to the KL loss {\rm (\ref{eqn:loss})} if
\begin{equation}\label{eqn:upper0}
a>-q+2,\quad b>2,\quad a+b \leq (q-r-1)/(2-v_w/v_0)-2r+2.
\end{equation}
There exist constants $a$ and $b$ satisfying (\ref{eqn:upper0}) if
$
q-r-1+(2-v_w/v_0)(q-2r-2)>0.
$
\end{prp}

Recall that $\ph_U(Y\mid X)$ is minimax and has a constant risk relative to the KL loss (\ref{eqn:loss}).
Hence if $\ph_{GB}(Y\mid X)$ is proper Bayes then it is admissible.
\begin{cor}
Assume that $q-r-1>0$.
Then $\ph_{GB}(Y\mid X)$ is admissible and minimax relative to the KL loss {\rm (\ref{eqn:loss})} if
\begin{equation}\label{eqn:upper1}
a>0,\quad  b>2, \quad a+b \leq (q-r-1)/(2-v_w/v_0)-2r+2.
\end{equation}
Thus, there exist constants $a$ and $b$ satisfying (\ref{eqn:upper1}) if $(-q+5r+1)/(2r)<v_w/v_0 <1$.
\end{cor}

Since $0<v_w<v_x\leq v_0$, it is observed that
$$
(q-r-1)/(2-v_w/v_0)-2r+2>(q-r-1)/2-2r+2=(q-5r+3)/2.
$$
We also obtain the following corollary.
\begin{cor}\label{cor:proper2}
Assume that $q-5r-1>0$.
Then, for any $v_x$, $v_y$ and $v_0\ (\geq v_x)$, $\ph_{GB}(Y\mid X)$ is admissible and minimax relative to the KL loss {\rm (\ref{eqn:loss})} if
$$
a>0,\quad  b>2, \quad a+b \leq (q-5r+3)/2.
$$
\end{cor}

\begin{figure}[tb]
\begin{center}
\includegraphics[bb=130 580 470 750]{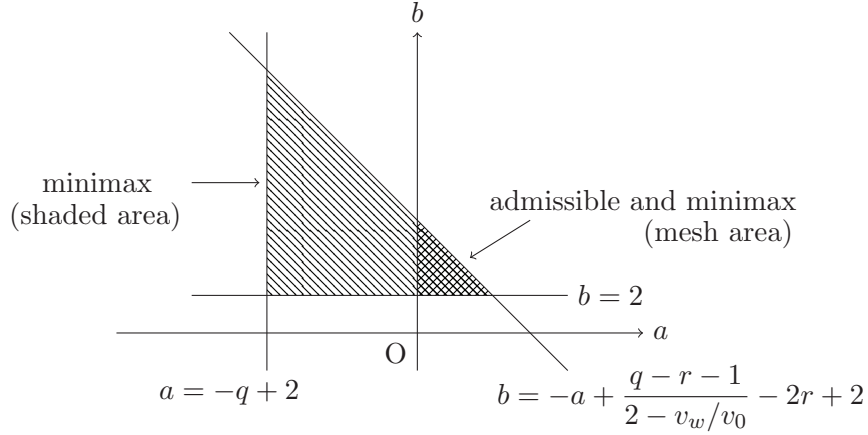}
\end{center}
\vspace{-15pt}
\caption{Sufficient conditions on $(a,b)$ of $\ph_{GB}(Y|X)$ for admissibility and minimaxity}
\end{figure}

{\bf Proof of Proposition \ref{prp:mini_GB}.}\ \ 
Using Theorem \ref{thm:faith}, we will derive a sufficient condition for minimaxity of $\ph_{GB}(Y\mid X)$.

\medskip
Denote
$$
c=a+b+2r-2.
$$
Let
\begin{align*}
\pi_{GB}^J(\La)&=v_1^{r(r+1)/2}|v_1I_r+(1-v_1)\La|^{-r-1}\pi_{GB}[\La\{v_1I_r+(1-v_1)\La\}^{-1}]\\
&=K_{a,b}v_1^{r(r+b-1)/2}|\La|^{a/2-1}|I_r-\La|^{b/2-1}|v_1I_r+(1-v_1)\La|^{-c/2}.
\end{align*}
When 
\begin{equation}\label{eqn:bound-cond}
q+a>2 \quad\textup{and}\quad b>2,
\end{equation}
it follows that for any $\La\in\partial\Rc_r$
$$
f_{GB}(\La;W)=(2\pi v)^{-qr/2}|\La|^{q/2}\pi_{GB}^J(\La)\exp\Big[-\frac{1}{2v}\tr(\La WW^\top)\Big]=0.
$$

\medskip
Define
\begin{align*}
m_{GB}&=m_{GB}(W)=\int_{\Rc_r} f_{GB}(\La;W)\dd\La, \\
M_{GB}&=M_{GB}(W)=\int_{\Rc_r}\La f_{GB}(\La;W)\dd\La.
\end{align*}
Since $0<v_1\leq 1$ and $\La\in\Rc_r$, it holds that
$$
|v_1I_r+(1-v_1)\La|^{-c/2}\leq \max\big(1,v_1^{-rc/2}\big),
$$
which implies that
$$
f_{GB}(\La;W)\leq \textup{const.}\times |\La|^{(q+a)/2-1}|I_r-\La|^{b/2-1}.
$$
Thus if $q+a>0$ and $b>0$ then $m_{GB}$ and $M_{GB}$ are finite.

\medskip
It is seen from Lemmas \ref{lem:diff1} and \ref{lem:diff3} that
$$
\frac{\Dc_\La \pi_{GB}^J(\La)}{\pi_{GB}^J(\La)}
=\frac{1}{2}\big[(a-2)\La^{-1}-(b-2)(I_r-\La)^{-1}-(1-v_1)c\{v_1I_r+(1-v_1)\La\}^{-1}\Big],
$$
so that
\begin{align*}
\De_1(W;\pi_{GB}^J)
&=4\int_{\Rc_r}\frac{1}{\pi_{GB}^J(\La)}\tr[\La^2\Dc_\La \pi_{GB}^J(\La)]f_{GB}(\La;W)\dd\La\\
&=2(a-2)\tr M_{GB}-2(b-2)\int_{\Rc_r}\tr[(I_r-\La)^{-1}\La^2]f_{GB}(\La;W)\dd\La\\
&\quad -2(1-v_1)c\int_{\Rc_r}\tr[\{v_1I_r+(1-v_1)\La\}^{-1}\La^2]f_{GB}(\La;W)\dd\La.
\end{align*}
Note that
$$
\tr[(I_r-\La)^{-1}\La^2]=-\tr\La+\tr[(I_r-\La)^{-1}\La],
$$
which leads to
\begin{align}\label{eqn:De1}
&\De_1(W;\pi_{GB}^J) \non\\
&=2(a+b-4)\tr M_{GB}-2(b-2)\int_{\Rc_r}\tr[(I_r-\La)^{-1}\La]f_{GB}(\La;W)\dd\La \non\\
&\quad -2(1-v_1)c\int_{\Rc_r}\tr[\{v_1I_r+(1-v_1)\La\}^{-1}\La^2]f_{GB}(\La;W)\dd\La.
\end{align}
Similarly, Lemmas \ref{lem:diff1} and \ref{lem:diff3} are used to see that
\begin{align*}
\frac{2\tr[\La M_{GB}\Dc_\La \pi_{GB}^J(\La)]}{\pi_{GB}^J(\La)}
&=(a-2)\tr M_{GB}-(b-2)\tr[M_{GB}(I_r-\La)^{-1}\La] \\
&\quad -(1-v_1)c\tr[M_{GB}\{v_1I_r+(1-v_1)\La\}^{-1}\La],
\end{align*}
which yields that
\begin{align}\label{eqn:De2}
\De_2(W;\pi_{GB}^J)
&= \frac{2}{m_{GB}}\int_{\Rc_r}\frac{\tr[\La M_{GB}\Dc_\La \pi_{GB}^J(\La)]}{\pi_{GB}^J(\La)}f_{GB}(\La;W)\dd\La \non\\
&= (a-2)\tr M_{GB} -\frac{b-2}{m_{GB}}\int_{\Rc_r}\tr[M_{GB}(I_r-\La)^{-1}\La]f_{GB}(\La;W)\dd\La \non\\
&\quad -\frac{(1-v_1)c}{m_{GB}}\int_{\Rc_r}\tr[M_{GB}\{v_1I_r+(1-v_1)\La\}^{-1}\La]f_{GB}(\La;W)\dd\La.
\end{align}
Hence combining (\ref{eqn:De1}) and (\ref{eqn:De2}) gives that
\begin{align}\label{eqn:De0}
\De(W;\pi_{GB}^J)&=\De_1(W;\pi_{GB}^J)-\De_2(W;\pi_{GB}^J)-(q-3r-3)\tr M_{GB} \non\\
&= -(q-3r+3-a-2b)\tr M_{GB} +\De_3+\De_4,
\end{align}
where
\begin{align*}
\De_3 &= -2(b-2)\int_{\Rc_r}\tr[(I_r-\La)^{-1}\La]f_{GB}(\La;W)\dd\La \\
&\qquad +\frac{b-2}{m_{GB}}\int_{\Rc_r}\tr[M_{GB}(I_r-\La)^{-1}\La]f_{GB}(\La;W)\dd\La, \\
\De_4 &=-2(1-v_1)c\int_{\Rc_r}\tr[\{v_1I_r+(1-v_1)\La\}^{-1}\La^2]f_{GB}(\La;W)\dd\La \\
&\qquad +\frac{(1-v_1)c}{m_{GB}}\int_{\Rc_r}\tr[M_{GB}\{v_1I_r+(1-v_1)\La\}^{-1}\La]f_{GB}(\La;W)\dd\La.
\end{align*}
Here, it can easily be verified that $\De(W;\pi_{GB}^J)$ is finite for $q+a>0$ and $b>2$.

\medskip
For notational simplicity, we use the notation 
$$
\Er_\La[g(\La)]=\int_{\Rc_r}g(\La)f_{GB}(\La;W)\dd\La \Big/\int_{\Rc_r}f_{GB}(\La;W)\dd\La
$$ 
for an integrable function $g(\La)$.
Then from (\ref{eqn:De0}),
\begin{align}\label{eqn:De00}
{\De(W;\pi_{GB}^J)\over m_{GB}}
=& (c-q+r+b-1)\tr \Er_\La(\La) \non\\
&+ (b-2)\Big[ \tr\big[ \Er_\La(\La)\Er_\La\{(I_r-\La)^{-1}\La\} \big] - 2\tr\big[ \Er_\La\{(I_r-\La)^{-1}\La\}\big]\Big]\non\\
&+(1-v_1)c\Big[ \tr\big[ \Er_\La(\La)\Er_\La\{(v_1I_r+(1-v_1)\La)^{-1}\La\}\big]
\non\\
&\qquad\qquad\qquad -2\tr\big[ \Er_\La\{(v_1I_r+(1-v_1)\La)^{-1}\La^2\}\big]\Big]
\end{align}
for $c=a+b+2r-2$.
Note that $0_{r\times r} \preceq \La \preceq I_r$ and $I_r \preceq (I_r-\La)^{-1}$.
Since $\Er_\La(\La) \preceq I_r$ and $\tr\big[ (I_r-\La)^{-1}\La \big] \geq \tr\La$, the second term in the r.h.s. of (\ref{eqn:De00}) is evaluated as
\begin{align*}
\tr\big[ &\Er_\La(\La)\Er_\La\{(I_r-\La)^{-1}\La\}\big] - 2\tr\big[ \Er_\La\{(I_r-\La)^{-1}\La\}\big]
\\
&\leq - \tr\big[ \Er_\La\{(I_r-\La)^{-1}\La\}\big] \leq - \tr \Er_\La(\La).
\end{align*}
Since $b>2$, we have
\begin{align}\label{eqn:De01}
{\De(W;\pi_{GB}^J)\over m_{GB}}
\leq & (c-q+r+1)\tr \Er_\La(\La) \non\\
&+(1-v_1)c\Big[ \tr\big[ \Er_\La(\La)\Er_\La\{(v_1I_r+(1-v_1)\La)^{-1}\La\}\big]
\non\\
&\qquad\qquad\qquad -2\tr\big[ \Er_\La\{(v_1I_r+(1-v_1)\La)^{-1}\La^2\}\big]\Big].
\end{align}
It is here observed that
\begin{align*}
&(1-v_1)\{v_1I_r+(1-v_1)\La\}^{-1}\La \\
&\qquad = \{v_1I_r+(1-v_1)\La\}^{-1} \{ v_1I_r+(1-v_1)\La - v_1I_r\}\\
&\qquad= I_r - v_1 \{v_1I_r+(1-v_1)\La\}^{-1},
\\
&(1-v_1)\{v_1I_r+(1-v_1)\La\}^{-1}\La^2 \\
&\qquad= \La - v_1 \{v_1I_r+(1-v_1)\La\}^{-1}\La,
\end{align*}
which is used to get
\begin{align*}
&(1-v_1)c\Big[ \tr\big[ \Er_\La(\La)\Er_\La\{(v_1I_r+(1-v_1)\La)^{-1}\La\}\big] \\
&\qquad\qquad -2\tr\big[ \Er_\La\{(v_1I_r+(1-v_1)\La)^{-1}\La^2\}\big]\Big] \\
&\qquad =-c\tr \Er_\La(\La)\\
&\qquad\qquad + cv_1 \Big[ 2 \tr\big[ \Er_\La\{(v_1I_r+(1-v_1)\La)^{-1}\La\}\big] \\
&\qquad\qquad\qquad-\tr\big[ \Er_\La(\La) \Er_\La\{(v_1I_r+(1-v_1)\La)^{-1}\}\big]\Big].
\end{align*}
Substituting this quantity into (\ref{eqn:De01}) gives
\begin{align}\label{eqn:De02}
{\De(W;\pi_{GB}^J)\over m_{GB}}
\leq & -(q-r-1)\tr \Er_\La(\La) \non\\
&+ cv_1 \Big[ 2 \tr\big[ \Er_\La\{(v_1I_r+(1-v_1)\La)^{-1}\La\}\big] \non\\
&\qquad\qquad -\tr\big[ \Er_\La(\La) \Er_\La\{(v_1I_r+(1-v_1)\La)^{-1}\}\big]\Big].
\end{align}

\medskip
To evaluate the second term in the r.h.s. of (\ref{eqn:De02}), note that
\begin{equation}
I_r \preceq \{v_1I_r+(1-v_1)\La\}^{-1} \preceq v_1^{-1}I_r.
\label{eqn:inq}
\end{equation}
In the case of $c\geq 0$, it is seen from (\ref{eqn:inq}) that
\begin{align*}
cv_1 &\Big[ 2 \tr\big[ \Er_\La\{(v_1I_r+(1-v_1)\La)^{-1}\La\}\big]-\tr\big[ \Er_\La(\La) \Er_\La\{(v_1I_r+(1-v_1)\La)^{-1}\}\big]\Big]
\\
&\leq cv_1 \Big\{ {2\over v_1}\tr \Er_\La(\La) - \tr \Er_\La(\La)\Big\}
= c (2-v_1) \tr \Er_\La(\La),
\end{align*}
which implies that
$$
\De(W;\pi_{GB}^J)/m_{GB}
\leq \{ - (q-r-1) + c(2-v_w/v_0)\}\tr \Er_\La(\La),
$$
because $1>v_1=v/v_0\geq v_w/v_0>0$.
It is noted that $c=a+b+2r-2>-q + 2r +2=- (q-r-1) + r+1$ because $a>2-q$ and $b>2$.
Thus, one gets a sufficient condition given by
\begin{equation}
\max \{0, - (q-r-1) + r+1\} \leq c \leq (q-r-1)/(2-v_w/v_0).
\label{eqn:sc1}
\end{equation}
In the case of $c\leq 0$,  it is seen from (\ref{eqn:inq}) that
\begin{align*}
cv_1 &\Big[ 2 \tr\big[ \Er_\La\{(v_1I_r+(1-v_1)\La)^{-1}\La\}\big]-\tr\big[ \Er_\La(\La) \Er_\La(v_1I_r+(1-v_1)\La)^{-1}\big]\Big]
\\
&\leq cv_1 \Big\{ 2\tr \Er_\La(\La) - {1\over v_1}\tr \Er_\La(\La)\Big\}
= c (2v_1-1) \tr \Er_\La(\La),
\end{align*}
which implies that
$$
\De(W;\pi_{GB}^J)/m_{GB}
\leq \{ - (q-r-1) + c(2v_w/v_0-1)\}\tr \Er_\La(\La).
$$
Hence, it holds true that $ - (q-r-1) + c(2v_w/v_0-1)\leq 0$ if 
\begin{equation}
\min \{0, - (q-r-1) + r+1\} \leq c \leq 0.
\label{eqn:sc2}
\end{equation}
Combining (\ref{eqn:sc1}) and (\ref{eqn:sc2}) yields the condition $- (q-r-1) + r+1 \leq c \leq (q-r-1)/(2-v_w/v_0)$, namely, $-q+4\leq a+b\leq (q-r-1)/(2-v_w/v_0)-2r+2$.
From (\ref{eqn:bound-cond}), the sufficient conditions on $(a,b)$ for minimaxity can be written as $a>2-q$, $b>2$ and $a+b\leq (q-r-1)/(2-v_w/v_0)-2r+2$ if
\begin{align*}
&\{(q-r-1)/(2-v_w/v_0)-2r+2\}-\{-q+4\}\\
&=\{q-r-1+(2-v_w/v_0)(q-2r-2)\}/(2-v_w/v_0)>0.
\end{align*}
Thus the proof is complete.
\hfill$\Box$

\begin{remark}
Take $v_x=v_w=v_0=1$.
Let $X|\Th\sim\Nc_{r\times q}(\Th,I_r\otimes I_q)$.
Consider the problem of estimating the mean matrix $\Th$ under the squared Frobenius norm loss $\Vert\Thh-\Th\Vert^2$.
The Bayesian estimator with respect to (\ref{eqn:pr_Th}) and (\ref{eqn:pr_GB}) is expressed as
$$
\Thh_{GB}=\bigg[I_r-\frac{\int_{\Rc_r}\Om|\Om|^{(q+a)/2-1}|I_r-\Om|^{b/2-1}\exp[-\tr(\Om XX^\top)/2]\dd\Om}{\int_{\Rc_r}|\Om|^{(q+a)/2-1}|I_r-\Om|^{b/2-1}\exp[-\tr(\Om XX^\top)/2]\dd\Om}\bigg]X.
$$
Then the same arguments as in this section yield that $\Thh_{GB}$ is proper Bayes and minimax if $a>0$, $b>2$, $q>3r+1$ and $ 2 < a+b \leq q-3r+1$.
\hfill$\Box$
\end{remark}

\section{Superharmonic priors for minimaxity}
\label{sec:superharmonic}

In estimation of the normal mean vector, Stein (1973, 1981) discovered an interesting relationship between superharmonicity of prior density and minimaxity of the resulting generalized Bayes estimator.
The relationship is very important and useful in Bayesian predictive density estimation.
In this section we derive some Bayesian minimax predictive densities with superharmonic priors.

\medskip
Let $\ph_\pi=\ph_\pi(Y|X)$ be a Bayesian predictive density with respect to a prior $\pi(\Th)$, where $\pi(\Th)$ is twice differentiable and the marginal density $m_\pi(X;v_x)$ is finite.
All the results in this section are based on the following key lemma.
\begin{lem}\label{lem:superharmonic}
Denote by $\nabla_\Th=(\partial/\partial\th_{ij})$ the $r\times q$ differentiation operator matrix with respect to $\Th$.
Then $\ph_\pi$ is minimax relative to the KL loss $(\ref{eqn:loss})$ if $\pi(\Th)$ is superharmonic, namely,
$$
\tr[\nabla_\Th\nabla_\Th^\top \pi(\Th)]=\sum_{i=1}^r\sum_{j=1}^q\frac{\partial^2 \pi(\Th)}{\partial \th_{ij}^2}\leq 0.
$$
\end{lem}

{\bf Proof.}\ \ 
This lemma can be proved along the same arguments as in Stein (1981).
See also George et al. (2006) and Brown et al. (2008).
\hfill$\Box$

\medskip
Define a class of prior densities as
$$
\pi(\Th)=g(\Si),\quad \Si=\Th\Th^\top,
$$
where $g$ is twice differentiable with respect to $\Si$.
Let $\Dc_\Si$ be an $r\times r$ matrix of differential operator with respect to $\Si=(\si_{ij})$ such that the $(i,j)$ element of $\Dc_\Si$ is
$$
\{\Dc_\Si\}_{ij}=\frac{1+\de_{ij}}{2}\frac{\partial}{\partial \si_{ij}},
$$
where $\de_{ij}$ stands for the Kronecker delta.
Let 
$$
G=(g_{ij})=G(\Si)=\Dc_\Si g(\Si),
$$
namely, $G$ is an $r\times r$ symmetric matrix such that $g_{ij}=\{\Dc_\Si\}_{ij} g(\Si)$.

\begin{lem}\label{lem:condition1}
$\ph_\pi$ with respect to $\pi(\Th)=g(\Si)$ is minimax relative to the KL loss $(\ref{eqn:loss})$ if
$$
\tr[\nabla_\Th\nabla_\Th^\top\pi(\Th)]=2[(q-r-1)\tr(G)+2\tr(\Dc_\Si \Si G)]\leq 0,
$$
where $G=\Dc_\Si g(\Si)$.
\end{lem}

{\bf Proof.}\ \ 
Using (i) and (ii) in Lemma \ref{lem:diff4} gives that
\begin{align*}
\tr[\nabla_\Th\nabla_\Th^\top \pi(\Th)]
&=2\tr(\nabla_\Th \Th^\top \Dc_\Si g(\Si))=2\tr(\nabla_\Th \Th^\top G)\\
&=2\big[(q-r-1)\tr(G)+2\tr(\Dc_\Si \Si G)\big].
\end{align*}
From Lemma \ref{lem:superharmonic}, the proof is complete.
\qed

\medskip
Let $\la_1,\ldots,\la_r$ be ordered eigenvalues of $\Si=\Th\Th^\top$, where $\la_1\geq\cdots\geq\la_r$, and let $\La=\diag(\la_1,\ldots,\la_r)$.
Denote by $\Ga=(\ga_{ij})$ an $r\times r$ orthogonal matrix such that $\Ga^\top\Si \Ga=\La$.
Assume that $g(\Si)$ is orthogonally invariant, namely, $g(\Si)=g(P\Si P^\top)$ for any orthogonal matrix $P$.
Then, we can assume that $g(\Si)=g(\La)$ without loss of generality.
\begin{prp}\label{prp:condition2}
Assume that $g(\Si)=g(\La)$ and $g(\La)$ is a twice differentiable function of $\La$.
Then $\ph_\pi$ with $\pi(\Th)=g(\La)$ is minimax relative to the KL loss $(\ref{eqn:loss})$ if
\begin{align*}
&\tr[\nabla_\Th\nabla_\Th^\top\pi(\Th)]\\
&=2\sum_{i=1}^r\bigg\{(q-r+1)\phi_i(\La)+\sum_{j\ne i}^r\frac{\la_i\phi_i(\La)-\la_j\phi_j(\La)}{\la_i-\la_j}+2\la_i\frac{\partial\phi_i(\La)}{\partial\la_i}\bigg\}\leq 0,
\end{align*}
where $\phi_i(\La)=\partial g(\La)/\partial\la_i$.
\end{prp}

{\bf Proof.}\ \ 
Since from (i) of Lemma \ref{lem:diff2}
$$
\{\Dc_\Si\}_{ij}\la_k=\ga_{ik}\ga_{jk},
$$
it is observed that by the chain rule
$$
\{\Dc_\Si\}_{ij} g(\La)=\sum_{k=1}^r \frac{\partial g(\La)}{\partial\la_k}\{\Dc_\Si\}_{ij}\la_k
=\{\Ga\Phi(\La)\Ga^\top\}_{ij},
$$
where $\Phi(\La)=\diag(\phi_1(\La),\ldots,\phi_r(\La))$.
Using Lemma \ref{lem:condition1} and (ii) of Lemma \ref{lem:diff2} gives that
\begin{align*}
&\tr[\nabla_\Th\nabla_\Th^\top\pi(\Th)] \\
&=2[(q-r-1)\tr\{\Ga\Phi(\La)\Ga^\top\}+2\tr\{\Dc_\Si \Ga\La\Phi(\La)\Ga^\top\}] \\
&=2\sum_{i=1}^r\bigg[(q-r-1)\phi_i(\La)+\sum_{j\ne i}^r\frac{\la_i\phi_i(\La)-\la_j\phi_j(\La)}{\la_i-\la_j}+2\frac{\partial}{\partial\la_i}\{\la_i\phi_i(\La)\}\bigg]\\
&=2\sum_{i=1}^r\bigg\{(q-r+1)\phi_i(\La)+\sum_{j\ne i}^r\frac{\la_i\phi_i(\La)-\la_j\phi_j(\La)}{\la_i-\la_j}+2\la_i\frac{\partial\phi_i(\La)}{\partial\la_i}\bigg\}.
\end{align*}
Hence the proof is complete.
\qed

\medskip
Using Proposition \ref{prp:condition2}, we give some examples of Bayesian predictive densities with respect to superharmonic priors.
Consider a class of shrinkage prior densities,
$$
\pi_{SH}(\Th) = \{\tr(\Th\Th^\top)\}^{-\be/2}\prod_{i=1}^r \la_i^{-\al_i/2} 
=\bigg\{\sum_{i=1}^r \la_i\bigg\}^{-\be/2}\prod_{i=1}^r \la_i^{-\al_i/2},
$$
where $\al_1,\ldots,\al_r$ and $\be$ are nonnegative constants.
The class $\pi_{SH}(\Th)$ includes both harmonic priors $\pi_{EM}(\Th)$ and $\pi_{JS}(\Th)$, which are given in (\ref{eqn:pr_em}) and (\ref{eqn:pr_js}), respectively.
Indeed, $\pi_{SH}(\Th)$ is the same as $\pi_{EM}(\Th)$ if $\al_1=\cdots=\al_r=\al^{EM}$ and $\be=0$ and as $\pi_{JS}(\Th)$ if $\al_1=\cdots=\al_r=0$ and $\be=\be^{JS}$.

\medskip
It is noted that
\begin{align}
\frac{\partial}{\partial \la_k}\pi_{SH}(\Th) &= -\frac{1}{2}\Big(\frac{\al_k}{\la_k}+\frac{\be}{\sum_{i=1}^r \la_i}\Big)\pi_{SH}(\Th), \label{eqn:d_pi_g} \\
\frac{\partial^2}{\partial \la_k^2}\pi_{SH}(\Th) &=\frac{1}{2}\bigg\{\Big(\frac{\al_k}{\la_k^2}+\frac{\be}{(\sum_{i=1}^r \la_i)^2}\Big)+\frac{1}{2}\Big(\frac{\al_k}{\la_k}+\frac{\be}{\sum_{i=1}^r \la_i}\Big)^2\bigg\}\pi_{SH}(\Th)  . \label{eqn:dd_pi_g}
\end{align}
Combining (\ref{eqn:d_pi_g}), (\ref{eqn:dd_pi_g}) and Proposition \ref{prp:condition2}, we obtain
\begin{align}
&\tr[\nabla_\Th\nabla_\Th^\top \pi_{SH}(\Th)] \non\\
&=\pi_{SH}(\Th)\sum_{i=1}^r\bigg[\{\al_i^2-(q-r-1)\al_i\}\frac{1}{\la_i}-2\sum_{j>i}^r\frac{\al_i-\al_j}{\la_i-\la_j}+\frac{2\al_i\be}{\tr(\Th\Th^\top)}\bigg] \non\\
&\qquad +\pi_{SH}(\Th)\frac{\be^2-(qr-2)\be}{\tr(\Th\Th^\top)}.
\label{eqn:dd-pi_g}
\end{align}

\begin{exm}\label{exm:1}
Let
$$
\pi_{ST}(\Th)=\prod_{i=1}^r \la_i^{-\al_i/2},
$$
where $\al_1,\ldots,\al_r$ are nonnegative constants.
Assume that $\al_1\geq\cdots\geq\al_r$.
Note that
\begin{align*}
\sum_{i=1}^r\sum_{j>i}^r\frac{\al_i-\al_j}{\la_i-\la_j}
&=\sum_{i=1}^r\sum_{j>i}^r\frac{1}{\la_i}\frac{\la_i-\la_j+\la_j}{\la_i-\la_j}(\al_i-\al_j)\\
&=\sum_{i=1}^r(r-i)\frac{\al_i}{\la_i}-\sum_{i=1}^r\frac{1}{\la_i}\sum_{j>i}^r\al_j+\sum_{i=1}^r\sum_{j>i}^r\frac{\la_j}{\la_i}\frac{\al_i-\al_j}{\la_i-\la_j}\\
&\geq \sum_{i=1}^r(r-i)\frac{\al_i}{\la_i}-\sum_{i=1}^r\frac{1}{\la_i}\sum_{j>i}^r\al_j.
\end{align*}
From (\ref{eqn:dd-pi_g}), it is seen that
\begin{align*}
&\tr[\nabla_\Th\nabla_\Th^\top \pi_{ST}(\Th)] \\
&\leq\pi_{ST}(\Th)\sum_{i=1}^r\bigg\{\al_i^2-(q-r-1)\al_i-2(r-i)\al_i+2\sum_{j>i}^r\al_j\bigg\}\frac{1}{\la_i}\\
&=\pi_{ST}(\Th)\sum_{i=1}^r\bigg\{\al_i^2-(q+r-2i-1)\al_i+2\sum_{j>i}^r\al_j\bigg\}\frac{1}{\la_i}.
\end{align*}
Here, assume additionally that $\al_i\leq\al_i^{ST}/2$ with $\al_i^{ST}=q+r-2i-1$ for $i=1,\ldots,r$.
For each $i$ we observe that
\begin{align*}
&\al_i^2-(q+r-2i-1)\al_i+2\sum_{j>i}^r\al_j \\
&\leq \al_{i+1}^2-(q+r-2i-1)\al_{i+1}+2\sum_{j>i}^r\al_j\\
&= \al_{i+1}^2-\{q+r-2(i+1)-1\}\al_{i+1}+2\sum_{j>i+1}^r\al_j\\
&\leq \cdots\\
&\leq \al_r^2-(q-r-1)\al_r\leq0,
\end{align*}
which implies that $\tr[\nabla_\Th\nabla_\Th^\top \pi_{ST}(\Th)]\leq 0$ if $\al_1\geq\cdots\geq\al_r$ and $\al_i\leq\al_i^{ST}/2$ for each $i$.
Then the resulting Bayesian predictive density is minimax under the KL loss (\ref{eqn:loss}).
\ \hfill$\Box$
\end{exm}

\begin{exm}\label{exm:2}
Consider a prior density of the form
\begin{equation}\label{eqn:pr_MS1}
\pi_{MS1}(\Th)=\{\tr(\Th\Th^\top)\}^{-\be^{MS}/2}\prod_{i=1}^r \la_i^{-\al_i^{ST}/4},
\end{equation}
where $\be^{MS}=2(r-1)$.
Combining Example \ref{exm:1} and (\ref{eqn:dd-pi_g}) gives that
\begin{align*}
&\tr[\nabla_\Th\nabla_\Th^\top \pi_{MS1}(\Th)]\\
&\leq\pi_{MS1}(\Th)\sum_{i=1}^r\frac{\al_i^{ST}\be^{MS}}{\tr(\Th\Th^\top)}
+\pi_{MS1}(\Th)\frac{(\be^{MS})^2-(qr-2)\be^{MS}}{\tr(\Th\Th^\top)}=0.
\end{align*}
Hence the Bayesian predictive density with respect to $\pi_{MS1}(\Th)$ is minimax relative to the KL loss (\ref{eqn:loss}).
\ \hfill$\Box$
\end{exm}

\medskip
In the literature, many shrinkage estimators have been developed in estimation of a normal mean matrix.
It is worth pointing out that the Bayesian predictive densities with superharmonic prior $\pi_{SH}(\Th)$ correspond to such shrinkage estimators.

\medskip
Let $X|\Th\sim\Nc_{r\times q}(\Th, v_x I_r\otimes I_q)$ and denote an estimator of $\Th$ by $\Thh$.
Consider the problem of estimating the mean matrix $\Th$ relative to quadratic loss $L_Q(\Thh,\Th)=\Vert\Thh-\Th\Vert^2$.
Then the generalized Bayes estimator of $\Th$ with the prior density $\pi_{SH}$ is expressed as
\begin{align*}
\Thh_{SH}
&=\frac{\int_{\Re^{r\times q}} \Th \exp(-\Vert X-\Th \Vert^2/(2v_x))\pi_{SH}(\Th)\dd\Th}{\int_{\Re^{r\times q}} \exp(-\Vert X-\Th \Vert^2/(2v_x))\pi_{SH}(\Th)\dd\Th}.
\end{align*}
If $\pi_{SH}$ is superharmonic then $\Thh_{SH}$ is minimax relative to the quadratic loss $L_Q$.

\medskip
Since $v_x\nabla_\Th \exp(-\Vert X-\Th \Vert^2/(2v_x))=-(\Th-X)\exp(-\Vert X-\Th \Vert^2/(2v_x))$, the integration by parts gives that
\begin{align*}
\Thh_{SH}
&=X-v_x\frac{\int_{\Re^{r\times q}} [\nabla_\Th \exp(-\Vert X-\Th \Vert^2/(2v_x))]\pi_{SH}(\Th)\dd\Th}{\int_{\Re^{r\times q}} \exp(-\Vert X-\Th \Vert^2/(2v_x))\pi_{SH}(\Th)\dd\Th} \\
&=X+v_x \frac{\int_{\Re^{r\times q}} \exp(-\Vert X-\Th \Vert^2/(2v_x))[\nabla_\Th\pi_{SH}(\Th)]\dd\Th}{\int_{\Re^{r\times q}} \exp(-\Vert X-\Th \Vert^2/(2v_x))\pi_{SH}(\Th)\dd\Th}.
\end{align*}
Here using (i) of Lemma \ref{lem:diff4} and (i) of Lemma \ref{lem:diff2} gives that
\begin{align*}
\nabla_\Th^\top \pi_{SH}(\Th)
&= 2\Th^\top \Dc_\Si \pi_{SH}(\Th) \non\\
&=- \Th^\top \bigg\{\Ga\diag\Big(\frac{\al_1}{\la_1},\ldots,\frac{\al_r}{\la_r}\Big)\Ga^\top+\frac{\be}{\tr(\Si)}I_r \bigg\} \pi_{SH}(\Th),
\end{align*}
which leads to
\begin{equation}\label{eqn:Th_MS}
\Thh_{SH}
=X-v_x\Er^{\Th|X}\bigg[ \bigg\{\Ga\diag\Big(\frac{\al_1}{\la_1},\ldots,\frac{\al_r}{\la_r}\Big)\Ga^\top+\frac{\be}{\tr(\Si)}I_r \bigg\} \Th\bigg],
\end{equation}
where $\Er^{\Th|X}$ stands for the posterior expectation with respect to a density proportional to $\exp(-\Vert \Th-X \Vert^2/(2v_x))\pi_{SH}(\Th)$.

\medskip
Denote by $XX^\top=HLH^\top$ the eigenvalue decomposition of $XX^\top$, where $H=(h_{ij})$ is an orthogonal matrix of order $r$ and $L=\diag(\ell_1,\ldots,\ell_r)$ is a diagonal matrix of order $r$ with $\ell_1\geq \cdots\geq \ell_r$.
Substituting $(X,H,L)$ for $(\Th,\Ga,\La)$ in the second term of the r.h.s. of (\ref{eqn:Th_MS}), we obtain an empirical Bayes shrinkage estimator
$$
\Thh_{MS}=X-v_x\bigg\{H\diag\Big(\frac{\al_1}{\ell_1},\ldots,\frac{\al_r}{\ell_r}\Big)H^\top+\frac{\be}{\tr(XX^\top)}I_r \bigg\}X.
$$

\medskip
The shrinkage estimator $\Thh_{MS}$ is equivalent to $\Thh_{JS}$, given in (\ref{eqn:JS}), when $\al_1=\cdots=\al_r=0$ and $\be=\be^{JS}$, and to $\Thh_{EM}$, given in (\ref{eqn:EM}), when $\al_1=\cdots=\al_r=\al^{EM}$ and $\be=0$.
In estimation of the normal mean matrix relative to the quadratic loss $L_Q$, $\Thh_{JS}$ and $\Thh_{EM}$ are minimax.

\medskip
If $\Thh_{MS}$ with certain specified $\al_1,\ldots,\al_r$ and $\be$ has good performance, the prior density $\pi_{SH}$ with the same $\al_1,\ldots,\al_r$ and $\be$ would produce a good Bayesian predictive density.
From Tsukuma (2008), $\Thh_{MS}$ is a minimax estimator dominating $\Thh_{EM}$  when $\al_i=\al_i^{ST}$ for $i=1,\ldots,r$ and $0\leq \be\leq 4(r-1)$.
A reasonable choice for $\be$ is $\be^{MS}=2(r-1)$ and this suggests that we should consider a prior density of the form
\begin{equation}\label{eqn:pr_MS2}
\pi_{MS2}(\Th)
=\{\tr(\Th\Th^\top)\}^{-\be^{MS}/2}\prod_{i=1}^r \la_i^{-\al_i^{ST}/2}
=\pi_{MS1}(\Th)\prod_{i=1}^r \la_i^{-\al_i^{ST}/4}.
\end{equation}
The prior density $\pi_{MS2}(\Th)$ is not superharmonic, and it is not known whether the resulting Bayesian predictive density is minimax or not.
In the next section, we verify risk behavior of the Bayesian predictive density with respect to $\pi_{MS2}(\Th)$ through Monte Carlo simulations.

\section{Monte Carlo studies}\label{sec:MCstudies}

This section briefly reports some numerical results so as to compare performance in risk of some Bayesian predictive densities for $r=2$ and $q=15$.

\medskip
First we investigate risk behavior of generalized Bayes predictive densities $\ph_{GB}(Y|X)$ with $v_0=1$ in the following six cases:
$$
(a,\, b)=(-11,\, 3),\ (-11,\, 9),\ (-11,\, 15),\ (-5,\, 3),\ (-5,\, 9),\ (1,\, 3)
$$
for the second-stage prior (\ref{eqn:pr_GB}).
When $r=2$ and $q=15$, $\ph_{GB}(Y|X)$ with the above six cases are minimax and, in particular, $\ph_{GB}(Y|X)$ with $(a,\, b)=(1,\,3)$ is proper Bayes for any $v_x$ and $v_y$ (see Corollary \ref{cor:proper2}).

\medskip
The risk has been simulated by 100,000 independent replications of $X$ and $Y$, where $X|\Th\sim\Nc_{r\times q}(\Th, v_xI_r\otimes I_q)$ and $Y|\Th\sim\Nc_{r\times q}(\Th, v_yI_r\otimes I_q)$ with $(v_x,v_y)=(0.1,\, 1),\ (1,\, 1)$ and $(1,\, 0.1)$.
It has been assumed that a pair of the maximum and the minimum eigenvalues of $\Th\Th^\top$ is $(0,\, 0),\ (24,\, 0)$ or $(24,\, 24)$.
Note that the best invariant predictive density $\ph_U(Y|X)$ has a constant risk and its risk is approximately given by
$$
R(\ph_U,\Th)=\frac{rq}{2}\log\frac{v_s}{v_y}\approx\begin{cases}
1.42 & \textup{for $(v_x,v_y)=(0.1,\ 1)$},\\
10.4 & \textup{for $(v_x,v_y)=(1,\ 1)$},\\
36.0 & \textup{for $(v_x,v_y)=(1,\ 0.1)$},
\end{cases}
$$
when $r=2$ and $q=15$.

\medskip
Denote by $\Bc(a,b)$ the matrix-variate beta distribution having the density (\ref{eqn:pr_GB}).
Using (\ref{eqn:m(W)-1}) with $\La=v_1\Om\{I_r-(1-v_1)\Om\}^{-1}$ and $v_1=v/v_0$, we can rewrite $\ph_{GB}(Y|X)$ as
$$
\ph_{GB}(Y|X)=\frac{\Er^{\Om}[g_{v_w}(\Om|W)]}{\Er^{\Om}[g_{v_x}(\Om|X)]}\ph_U(Y|X),
$$
where $\Er^{\Om}$ indicates expectation with respect to $\Om\sim \Bc(a+q,b)$ and 
$$
g_{v}(\Om|Z)=\Big|I_r-\Big(1-\frac{v}{v_0}\Big)\Om\Big|^{-q/2}\exp\Big[-\frac{1}{2v_0}\tr\Big[\Om\Big\{I_r-\Big(1-\frac{v}{v_0}\Big)\Om\Big\}^{-1}ZZ^\top\Big]\Big]
$$
for an $r\times q$ matrix $Z$.
Hence in our simulations, the expectation $\Er^{\Om}[g_{v}(\Om|Z)]$ was estimated by $j_0^{-1}\sum_{j=1}^{j_0} g_{v}(\Om_j|Z)$, where $j_0=100,000$ and the $\Om_j$ are independent replications from $\Bc(a+q,b)$.

\begin{table}[tb]
\begin{center}
{\scriptsize
\caption{Some simulated risk of generalized Bayes minimax predictive densities for $v_0=1$.}
\label{tab:1}
\vspace{2pt}
$
\begin{array}{ccccccccc}
\hline
(v_x, v_y) & {\rm Eigenvalues}&{\rm Minimax}&\multicolumn{6}{c}{(a,b)}\\
\cline{4-9}
& {\rm of}\ \Th\Th^\top &{\rm risk}&(-11,3)&(-11,9)&(-11,15)&(-5,3)&(-5,9)&(1,3) \\
\hline
(0.1,1) &(\ 0,\ 0)&1.42& 0.47 & 0.96 & 1.20 & 0.38 & 0.80 & 0.33 \\
        &(24,\ 0) &    & 0.91 & 1.17 & 1.30 & 0.87 & 1.09 & 0.84 \\
        &(24,24)  &    & 1.39 & 1.39 & 1.40 & 1.37 & 1.37 & 1.39 \\
[6pt]                   
(1, 1) &(\ 0,\ 0)&10.4&  5.3 &  6.9 &  7.7 &  2.8 &  4.6 &  1.6 \\
       &(24,\ 0) &    &  6.9 &  8.0 &  8.4 &  5.3 &  6.3 &  4.6 \\
       &(24,24)  &    &  8.7 &  9.0 &  9.2 &  7.9 &  8.0 &  8.4 \\
[6pt]                   
(1,0.1)&(\ 0,\ 0)&36.0& 15.2 & 24.0 & 28.2 &  9.6 & 17.7 &  6.8 \\
       &(24,\ 0) &    & 23.6 & 28.5 & 30.9 & 20.1 & 24.5 & 18.7 \\
       &(24,24)  &    & 32.7 & 33.2 & 33.5 & 31.0 & 31.4 & 32.4 \\
\hline
\end{array}
$
}
\end{center}
\end{table}

\medskip
The simulated results for risk of $\ph_{GB}(Y|X)$ are given in Table \ref{tab:1}.
When the pair of eigenvalues of $\Th\Th^\top$ is $(0,\, 0)$, our simulations suggest that the risk of $\ph_{GB}(Y|X)$ decreases as $a$ increases under which $b$ is fixed or under which $a+b$ is fixed and also that the risk of $\ph_{GB}(Y|X)$ increases as $b$ increases under which $a$ is fixed.
It is observed that $\ph_{GB}(Y|X)$ with $(a, b)=(1, 3)$ is superior to others.

\medskip
When the pair of eigenvalues of $\Th\Th^\top$ is $(24,\, 24)$, $\ph_{GB}(Y|X)$ with $(a, b)=(-5, 3)$ or $(-5, 9)$ is best, but the improvement over $\ph_U(Y|X)$ is little.
When the pair of eigenvalues of $\Th\Th^\top$ is $(24,\, 0)$, $\ph_{GB}(Y|X)$ with $(a, b)=(1, 3)$ is best.

\medskip
Next, we investigate the risk of Bayesian predictive densities based on superharmonic priors when $r=2$ and $q=15$.
If $\pi_s(\Th)$ is a superharmonic prior, then the Bayesian predictive density (\ref{eqn:BPD}) can be expressed as
$$
\ph_{\pi_s}(Y|X)=\frac{\Er^{\Th|W}[\pi_s(\Th)]}{\Er^{\Th|X}[\pi_s(\Th)]}\ph_U(Y|X),
$$
where $\Er^{\Th|W}$ and $\Er^{\Th|X}$ stand, respectively, for expectations with respect to $\Th|W\sim\Nc_{r\times q}(W, v_w I_r\otimes I_q)$ and $\Th|X \sim\Nc_{r\times q}(X, v_x I_r\otimes I_q)$.
In our simulations, $\ph_{\pi_s}(Y|X)$ was estimated by means of
$$
\ph_{\pi_s}(Y|X)\approx \frac{\sum_{i=1}^{i_0}\pi_s(\Th_i)}{\sum_{i=1}^{i_0}\pi_s(\Th_i)}\ph_U(Y|X),
$$
where $i_0=100,000$ and the $\Th_i$ and the $\Th_j$ are, respectively, independent replications from $\Nc_{r\times q}(W, v_w I_r\otimes I_q)$ and $\Nc_{r\times q}(X, v_x I_r\otimes I_q)$.

\begin{table}[bt]
\begin{center}
{\scriptsize
\caption{Some simulated risk of Bayesian predictive densities.}
\label{tab:2}
\vspace{2pt}
$
\begin{array}{cccc@{\hspace{20pt}}ccccc}
\hline
v_x & v_y & {\rm Eigenvalues} &{\rm Minimax}& $GB$ & $JS$\ & $EM$\, & $MS1$ & $MS2$ \\
&& {\rm of}\ \Th\Th^\top &{\rm risk}&&&&&\\
\hline
0.1& 1 &(\ 0,\ 0)&1.42& 0.33 & 0.09 & 0.28 & 0.71 & 0.09 \\
   &   &(24,\ 0) &    & 0.84 & 1.30 & 0.83 & 1.11 & 0.82 \\
   &   &(24,\ 4) &    & 1.27 & 1.31 & 1.27 & 1.30 & 1.26 \\
   &   &(24,\ 8) &    & 1.32 & 1.33 & 1.33 & 1.34 & 1.32 \\
   &   &(24, 12) &    & 1.35 & 1.34 & 1.35 & 1.36 & 1.35 \\
   &   &(24, 24) &    & 1.39 & 1.36 & 1.37 & 1.38 & 1.37 \\
[6pt]         
 1 & 1 &(\ 0,\ 0)&10.4&  1.6 &  0.7 &  2.1 &  5.2 &  0.7 \\
   &   &(24,\ 0) &    &  4.6 &  5.5 &  4.9 &  7.0 &  4.4 \\
   &   &(24,\ 4) &    &  5.7 &  5.9 &  5.9 &  7.4 &  5.4 \\
   &   &(24,\ 8) &    &  6.6 &  6.3 &  6.6 &  7.7 &  6.2 \\
   &   &(24, 12) &    &  7.2 &  6.6 &  7.1 &  8.0 &  6.7 \\
   &   &(24, 24) &    &  8.4 &  7.3 &  7.9 &  8.4 &  7.6 \\
[6pt]                    
 1 &0.1&(\ 0,\ 0)&36.0&  6.8 &  2.4 &  7.2 & 18.0 &  2.4 \\
   &   &(24,\ 0) &    & 18.7 & 25.8 & 18.9 & 26.1 & 18.2 \\
   &   &(24,\ 4) &    & 25.5 & 26.8 & 25.7 & 28.9 & 25.0 \\
   &   &(24,\ 8) &    & 28.1 & 27.7 & 28.2 & 30.2 & 27.6 \\
   &   &(24, 12) &    & 29.7 & 28.4 & 29.5 & 30.9 & 28.9 \\
   &   &(24, 24) &    & 32.4 & 30.0 & 31.3 & 32.0 & 30.8 \\
\hline
\end{array}
$
}
\end{center}
\end{table}

\medskip
The risk is based on 100,000 independent replications of $X$ and $Y$ for some pairs of two eigenvalues of $\Th\Th^\top$.
The simulation results are provided in Table \ref{tab:2}, where GB, JS, EM, MS1 and MS2 are the Bayesian predictive densities with the following priors.
\begin{list}{}{
\topsep=4pt
\parsep=0pt
\parskip=0pt
\itemsep=4pt
\itemindent=0pt
\labelwidth=40pt
\labelsep=8pt
\leftmargin=48pt
\listparindent=12pt
}
\item[GB:] (\ref{eqn:pr_Th}) and (\ref{eqn:pr_GB}) with $a=1$, $b=3$ and $v_0=1$,
\item[JS:] (\ref{eqn:pr_js}),
\item[EM:] (\ref{eqn:pr_em}),
\item[MS1:] (\ref{eqn:pr_MS1}),
\item[MS2:] (\ref{eqn:pr_MS2}).
\end{list}
Note that GB, JS, EM and MS1 are minimax, while MS2 has not been shown to be minimax.

\medskip
When the pair of eigenvalues of $\Th\Th^\top$ is $(0,\, 0)$, JS and MS2 are superior.
When the pair of eigenvalues of $\Th\Th^\top$ is $(24,\, 24)$, JS has nice performance but it is bad if the two eigenvalues of $\Th\Th^\top$ are much different.

\medskip
Our simulations suggest that MS2 is better than EM and MS1.
When the two eigenvalues of $\Th\Th^\top$ are much different, namely they are $(24,\, 0)$ and $(24,\, 4)$, MS2 is best and GB or EM is second-best.

\bigskip
\noindent
{\bf Acknowledgments.}\ \
The research of the first author was supported by Grant-in-Aid for Scientific Research (15K00055) from Japan Society for the Promotion of Science (JSPS).
The research of the second author was supported in part by Grant-in-Aid for Scientific Research (15H01943 and 26330036) from JSPS.

\end{document}